\documentclass[12pt]{article}
\usepackage{amsmath,amsthm,amsfonts,amscd,amssymb,eucal,bbm}
\usepackage[margin=3cm]{geometry}
\begin{document}
\def\RR{{\mathbb R}}
\def\CC{{\mathbb C}}
\def\NN{{\mathbb N}}
\def\ZZ{{\mathbb Z}}
\def\s1{{S^1}}
\def\diff{{{\rm Diff}^+(S^1)}}
\def\vect{{{\rm Vect}(S^1)}}
\def\sl2{{{\rm SL}(2,\RR)}}
\def\psl2{{{\rm PSL}(2,\RR)}}
\def\mob{\text {M\"{o}b}}
\def\u1{{{\rm U}(1)}}
\def\su2{{{\rm SU}(2)}}
\def\suN{{{\rm SU}(N)}}
\def\so3{{{\rm SO}(3)}}
\def\gk {{G_k}}
\def\Dl{{\mathcal D}(L_0)}
\newcommand{\fin}{{\mathcal D}\!_{f\!in}}
\def\A{{\mathcal A}}
\def\B{{\mathcal B}}
\def\C{{\mathcal C}}
\def\D{{\mathcal D}}
\def\F{{\mathcal F}}
\def\H{{\mathcal H}}
\def\I{{\mathcal I}}
\def\K{{\mathcal K}}
\def\k{{\rm K}}
\def\M{{\mathcal M}}
\def\N{{\mathcal N}}
\def\O{{\mathcal O}}
\def\P{{\mathcal P}}
\def\R{{\mathcal R}}
\def\S{{\mathcal S}}
\def\T{{\mathcal T}}
\def\U{{\mathcal U}}
\def\V{{\mathcal V}}
\def\W{{\mathcal W}}
\def\G{{\bf G}}

\newtheorem{theorem}{Theorem}[section]
\newtheorem{definition}[theorem]{Definition}
\newtheorem{corollary}[theorem]{Corollary}
\newtheorem{proposition}[theorem]{Proposition}
\newtheorem{lemma}[theorem]{Lemma}
\newtheorem{remark}[theorem]{Remark}
\title{{\huge {\bf On the Uniqueness of Diffeomorphism Symmetry 
in Conformal Field Theory}}}

\author{SEBASTIANO CARPI 
\footnote{Supported in part by the Italian MIUR and GNAMPA-INDAM.} 
\\ Dipartimento di Scienze 
\\ Universit\`a ``G. d'Annunzio" di Chieti-Pescara 
\\ Viale Pindaro 87, I-65127 Pescara, Italy
\\ E-mail: carpi@sci.unich.it 
\\ {} \\ MIH\'ALY WEINER 
\\ Dipartimento di Matematica 
\\ Universit\`a di Roma ``Tor Vergata'' 
\\ Via della Ricerca Scientifica 1, I-00133 Roma, Italy
\\ E-mail: weiner@mat.uniroma2.it}
\date{} 

\maketitle 
\renewcommand{\sectionmark}[1]{} 

\begin{abstract} 
A M\"obius covariant net of von Neumann algebras 
on $S^1$ is diffeomorphism covariant if 
its M\"obius symmetry extends to diffeomorphism symmetry. 
We prove that in case the net is either a Virasoro net or any at least
$4$-regular net such an extension is unique: the local algebras
together with the M\"obius symmetry (equivalently: the local 
algebras together with the vacuum vector) completely determine it. 
We draw the two following conclusions for such theories. 
(1) The value of the central charge $c$ is an invariant and hence
the Virasoro nets for different values of $c$ are not isomorphic as 
M\"obius covariant nets.
(2) A vacuum preserving internal symmetry always commutes 
with the diffeomorphism symmetries. 
We further use our result to give a large class of new
examples of nets (even strongly additive ones), which are not
diffeomorphism covariant; i.e. which do not admit an 
extension of the symmetry to $\diff$. 
\end{abstract}

\section{Introduction}

This paper is motivated by the following question: is the $\diff$ symmetry, 
or the corresponding Virasoro algebra symmetry, exhibited by
2-dimensional Quantum Field Theory models unique?

We shall give a precise formulation to this question in the framework of 
Algebraic Quantum Field Theory (see the book of R. Haag \cite{Haag}). 
In this framework a chiral 2-dimensional quantum field theory is commonly 
described by means of a M\"obius covariant net of von Neumann algebras 
on $\s1$. The net is said to be diffeomorphism covariant if the 
corresponding positive energy representation of the M\"obius group 
has an extension to a (strongly continuous) projective unitary
representation of $\diff$ that acts covariantly on the von Neumann
algebras associated to the intervals of $\s1$ and that is compatible 
with the local structure of the net (see Sect. \ref{sectprel} for the
precise definition). 
It is known that this extension does not exist in general 
(see e.g. \cite{GLW,koester03a} and cf. also Sect. \ref{secttensorproduct} 
below) but to the best of our knowledge no results about its uniqueness 
appears in the literature despite the fact that this problem seems to be 
very natural. 
Besides of its mathematical naturalness the relevance of the above
uniqueness is strengthened by the increasing importance played in the past 
years by diffeomorphism symmetry of nets of von Neumann algebras on $\s1$ 
in the investigation of the structural properties of two-dimensional conformal 
field theories, see e.g. \cite{KL,KL2,koester03b,xu 2003,Carpi03b,LX}.

The main result of this paper is the proof that for large class of
diffeomorphism covariant nets on $\s1$ the $\diff$ symmetry is unique
in the sense explained above and hence that it is completely determined by
the underlying structure of the M\"obius covariant net. More precisely we
prove uniqueness for all Virasoro nets, namely the nets generated by the
zero-energy representations of $\diff$, (Theorem \ref{uniqvir}) and for
all $4$-regular diffeomorphism covariant nets on $\s1$ (Theorem
\ref{uniq4-reg}).  The latter class (see Sect. \ref{sectprel} for the
definition) includes every strongly additive diffeomorphism covariant net
on $\s1$ and hence every diffeomorphism covariant net which is completely
rational in the sense of \cite{KLM}, the nets generated by chiral current
algebras \cite{BS-M,GLW,Tol97,Was A} and their orbifold subnets \cite{xu
2003}.  
Since the M\"obius symmetry of a given net on $\s1$ is completely
determined by the vacuum vector \cite[Theorem 2.19]{FrG} our result shows
that in the above cases the $\diff$ symmetry of the net is also determined
by this vector.  Note also that the known examples of M\"obius covariant
nets which are not $4$-regular are not diffeomorphism covariant (see
\cite{GLW,koester03a}), so our uniqueness result could apply to every
diffeomorphism covariant net on $\s1$.

Let us now discuss some consequences of our results. Firstly the
uniqueness in the case of Virasoro nets implies that two Virasoro nets
cannot be isomorphic as M\"obius covariant nets on the circle if they have
different central charges (Corollary \ref{isovirnets}), a fact that seems
to be widely expected (see e.g. the introduction of \cite{BS-M}) but
that has not been explicitly stated in the literature. Similarly two
$4$-regular diffeomorphism covariant nets cannot be isomorphic as
M\"obius covariant nets on $\s1$ if the corresponding representations of
$\diff$ are not unitarily equivalent and in particular if they have a
different central charge (Corollary \ref{invariants}).  Another
interesting consequence is that we have a model independent proof of the
fact that diffeomorphisms symmetries commute with vacuum preserving
internal symmetries of a given $4$-regular net (Corollary \ref{com_int}).
Finally we apply our result to show that the tensor product of an infinite
sequence of $4$-regular diffeomorphism covariant net on $\s1$ is not
diffeomorphism covariant (Theorem \ref{difftensor}).

This paper is organized as follows. In Sect. \ref{sectprel} we discuss
various preliminaries about M\"obius covariant nets on $\s1$, 
subnets and diffeomorphism covariance together with its relation to the
Virasoro algebra. Almost all this facts are already carefully discussed in
the literature and we include them only to fix the notation and to keep
the paper reasonably self-contained. In Sect. \ref{sectuniqvir} we prove
the uniqueness of the $\diff$ symmetry in the case of Virasoro nets. The
result is obtained by showing in a rather direct way that the (chiral)
stress-energy tensors associated to two representations of $\diff$ making
a Virasoro net diffeomorphism covariant have to coincide.  In Sect.
\ref{sectestimates} we show that the maps corresponding to projective
unitary representations of $\diff$ continuously extend to a certain family
of nonsmooth diffeomorphisms in an appropriate topology. The result is 
proved at the Lie algebra level. Since the estimates in the paper
of Goodmann and Wallach \cite{GoWa} are not sufficient for our purpose we
need a more detailed analysis. In particular we cannot use directly
Nelson's commutator theorem \cite{Ne,RSII} to show that the operators
involved are self-adjoint but we find an estimate involving
the contraction semigroup associated to the conformal Hamiltonian $L_0$
which allows us to demonstrate self-adjointness following the 
ideas of the paper of E. Nelson \cite{Ne}. In Sect. \ref{sectuniq4reg} we
use the results of Sect. \ref{sectestimates} to construct a nontrivial local
operator which belongs to the Virasoro subnet associated to an arbitrary
representation of $\diff$ making a given $4$-regular net diffeomorphism
covariant. This construction, together with the main result in Sect.
\ref{sectuniqvir} and the minimality property of Virasoro nets proved in
\cite{Carpi98} allows us to reach the main objective of this paper, 
namely the uniqueness of the $\diff$ symmetry for $4$-regular nets. 
Finally in Sect. \ref{secttensorproduct} we discuss the above mentioned 
application of our main result to the case of infinite tensor products of 
nets.

\section{Preliminaries}
\label{sectprel}
\subsection{M\"obius covariant nets}
Let $\I$ be the set of open, nonempty and nondense arcs 
(also called: open proper arcs or open proper intervals) 
of the unit circle $S^1 =\{z\in \CC : \,|z|=1 \}$.
A {\bf M\"obius covariant net on $S^1$} is a map $\A$
which assigns to every open proper arc 
$I \subset S^1$ a von Neumann algebra $\A(I)$ acting on a fixed 
complex, infinite dimensional separable Hilbert space $\H_\A$ 
(``the vacuum Hilbert space
of the theory''), together with a given strongly continuous 
representation $U$ of $\mob \simeq \psl2$, the group of M\"obius
transformations\footnote{diffeomorphisms of $\s1$ of the form 
$z \mapsto \frac{az+b}{\overline{b}z+\overline{a}}$ with $a,b\in \CC$, 
$|a|^2-|b|^2=1$.} of the unit circle $S^1$ satisfying for all $I_1,I_2,I
\in \I$ and $\varphi \in \mob$ the following properties: 

\begin{itemize}
\item[(i)] {\it Isotony.}
\begin{equation}
I_1 \subset I_2 \,\Rightarrow\,
\A(I_1) \subset \A (I_2),
\end{equation}

\item[(ii)] {\it Locality.}
\begin{equation}
I_1 \cap I_2 = \emptyset \,\Rightarrow\,
[\A(I_1),\A(I_2)]=0,
\end{equation}

\item[(iii)] {\it Covariance.} 
\begin{equation}
U(\varphi){\A}(I)U(\varphi)^{-1}={\A}(\varphi(I)),
\end{equation}
 
\item[(iv)] {\it Positivity of the energy.} 
The representation $U$ is of positive energy type: 
the conformal Hamiltonian $L_0$, defined by 
$U(\theta_\alpha)=e^{i\alpha L_0}$ where $\theta_\alpha
\in \mob$ is the anticlockwise rotation by degree $\alpha$,
is positive.

\item[(v)] {\it Existence and uniqueness of the vacuum.} 
There exists a unique (up to phase) unit vector 
$\Omega \in \H_\A$ called the ``vacuum vector''
which is invariant under the action of $U$.
(Equivalently: up to phase there exists a unique
unit vector $\Omega$ that is of {\it zero-energy} 
for $U$; i.e. eigenvector of $L_0$ with eigenvalue zero.)

\item[(vi)] {\it Cyclicity of the vacuum.} 
${\Omega}$ is cyclic for the algebra 
$\A(S^1)=\bigvee_{I\in \I} \A(I)$.

\end{itemize} 

Some consequences of the axioms are \cite{FrG,GuLo96,FJ}:
\begin{itemize}
\item[(i)] {\it Reeh-Schlieder property.}
$\Omega$ is a cyclic and separating vector of the algebra $\A(I)$
for every $I \in \I$.

\item[(ii)] {\it Bisognano-Wichmann property.} 
\begin{equation}
U(\Lambda_I ({2\pi t}))=\Delta_I^{it}
\end{equation}  
where ${\Delta}_{I}$ is the modular operator associated to 
$\B(I)$ and $\Omega$, and $\Lambda_I$ is the one-parameter group 
of M\"obius transformations preserving the interval 
$I$ (the dilations associated to $I$) with the ``right'' 
parametrization (see e.g. \cite{GLW}).

\item[(iii)] {\it Haag duality.} For every $I \in \I$
\begin{equation}
\A(I)'= \A(I'),
\end{equation}
where $I'$ denotes the interior of the complement set of 
$I$ in $S^1$.

\item[(iv)] {\it Irreducibility.} $\A(S^1) = \bigvee_{I \in \I} \A(I)
={\rm B}(\H_\A)$, where ${\rm B}(\H_\A)$ denotes the algebra of all
bounded linear operators on $\H_\A$.

\item[(v)]{\it Additivity.} 
If $\S\subset \I$ is a covering of the 
interval $I$ then  
\begin{equation}
\A(I)\subset \bigvee_{J\in \S}\A(J).
\end{equation}
\end{itemize}

As a consequence of the {\it Bisognano-Wichmann property}, since $\mob$ is
generated by the dilations (associated to different intervals), the
representation $U$ is completely determined by the local algebras and the
vacuum vector via modular structure. Thus, we may say that there is a kind
of uniqueness regarding the representation of the M\"obius group.

According to the last property (additivity) and the isotony, 
if $I_1,I_2,I \in \I$ are such that $I_1 \cup I_2 = I$ then
$\A(I_1)\vee\A(I_2)=\B(I)$. In many (but not all)
physically interesting model an even stronger 
additivity property holds. The net $\A$ is said to be {\bf strongly
additive}, if $\A(I_1)\vee\A(I_2)=\A(I)$ whenever $I_1,I_2$ are the
connected components of $I\setminus \{p\}$ where $p$ is a point of 
the open interval $I$.

For an $n=2,3,..$ the net $\A$ is said to be {\bf n-regular}, if 
whenever we remove $n$ points from the circle the algebras associated to
the remaining intervals generate the whole of $\A(S^1)={\rm B}(\H_\A)$. 
By isotony $n$-regularity is a stronger property then $m$-regularity
if $n>m$, and by Haag duality every M\"obius covariant net is at least 
$2$-regular. Strong additivity is of course stronger than 
$n$-regularity for any $n$. 

All these properties are indeed ``additional'': there
are M\"obius covariant nets which are not even $3$-regular (see the 
examples in \cite{GLW}).

\subsection{Diffeomorphism covariance and the Virasoro nets}
\label{subsec:diffcov}

Let $\diff$ be the group of orientation preserving (smooth) 
diffeomorphisms of the circle. It is an infinite dimensional Lie group
with Lie algebra
identified with the real topological vector space $\vect$ of smooth real
vectors fields on $S^1$ with the usual $C^\infty$ topology
\cite[Sect. 6]{Milnor} endowed with the bracket given by the negative of
the usual brackets of vector fields. In this paper often we shall think of
a the vector field symbolically written as
$f(e^{i\vartheta})\frac{d}{d\vartheta}\in \vect$ as the corresponding real
function $f$. Also we shall use the notation $f'$ (calling it simply the
derivative) for the function on the circle obtained by derivating with
respect to the angle:  
$f'(e^{i\theta})=\frac{d}{d\alpha}f(e^{i\alpha})|_{\alpha=\theta}$.

A strongly continuous projective unitary representation 
$V$ of $\diff$ on a Hilbert space $\H$ is a strongly 
continuous $\diff \rightarrow \U(\H)/\mathbb T$ homomorphism.
The restriction of $V$ to $\mob \subset \diff$ 
always lifts to a unique strongly continuous unitary
representation of the universal covering group 
$\widetilde{\mob}$ of $\mob$. $V$ is said to be of positive energy type,
if its conformal Hamiltonian $L_0$, defined by the above representation
of $\widetilde{\mob}$ (similarly as in case of a representation of the
group $\mob$) has nonnegative spectrum. 
In this case we shall simply say that $V$ is a positive energy 
representation of $\diff$.

Sometimes for a $\gamma \in \diff$ we shall think of $V(\gamma)$ as a
unitary operator. Although there are more than one ways to fix the phases,
note that expressions like Ad$(V(\gamma))$ or $V(\gamma) \in \M$ for a von
Neumann algebra $\M \subset {\rm B}(\H)$ are unambiguous. We shall also say
that $V$ is an extension of the unitary representation $U$ of $\mob$ if we
can arrange the phases in such a way that $V(\varphi)=U(\varphi)$, or
without mentioning phases: Ad$(V(\varphi))=$ Ad$(U(\varphi))$, for all
$\varphi \in \mob$.

We have to keep in mind that after choosing phases the equality
of $V$ to another projective representation
$\tilde{V}$ means that $V(\gamma)^*\tilde{V}(\gamma)$ is a
multiple of the identity (and not necessary {\it the} identity) 
for all $\gamma \in \diff$.

\begin{definition}
\label{diffcov:def}
A M\"obius covariant net $(\A,$U$)$ is {\bf diffeomorphism
covariant} if there is a strongly continuous projective 
unitary representation 
$V$ of $\diff$ on $\H_\A$ such that for all 
$\gamma \in \diff$ and $I,J \in \I$ 
\begin{itemize} 
\item[1.]
$\gamma \in \mob \Rightarrow $ Ad$(V(\gamma)) = $ Ad$(U(\gamma))$
\item[2.]
$\gamma|_I={\rm{id}}_I \Rightarrow
\rm{Ad}(V(\gamma))|_{\A(I)}=\rm{id}_{\A(I)}$.
\item[3.]
$\gamma(I)=J \Rightarrow V(\gamma)\A(I)V(\gamma)^* = \A(J).$
\end{itemize}  
In particular $V$ is a positive energy representation of $\diff$ extending 
$U$. 
\end{definition}

Note that as a consequence of {\it Haag duality} and of the second 
of the above listed properties, if a diffeomorphism localized in
the interval $I$ --- i.e. it acts trivially (identically) elsewhere ---
then the corresponding unitary is also localized in $I$ in the sense that
it belongs to $\A(I)$.
 
The majority of the known examples of ``interesting'' conformal field
theories are diffeomorphism covariant. In fact it may turn out to be that
under some ``regularity" condition imposed on the net diffeomorphism
covariance is automatic. In this respect the examples of
non-diffeomorphism covariant
nets which we shall give in the last section are useful in showing that
for example strong additivity in itself is not a sufficient condition. (Up
to the knowledge of the authors, there have been no previous examples of
strongly additive nets that are not diffeomorphism covariant.)

We now briefly describe the irreducible positive energy representations of
$\diff$ --- for fixing notations rather than to introduce them ---
and the so-called Virasoro nets. (Find more in \cite{KaRa},\cite{GoWa},
\cite{Tol99} and
\cite{DMS}, for example.) For certain values of the central charge $c>0$
and the lowest weight $h \geq 0$ there is a positive energy
projective representation denoted by $V_{(c,h)}$ on the
Hilbert space $\H_{(c,h)}$.
In $V_{(c,h)}$ the spectrum of the conformal Hamiltonian
Sp$(L_0)=\{h,h+1,h+2,..\}$ unless $h=0$ in which case  
the value $h+1=1$ is missing from it; all these corresponding
of course to eigenvalues, only. The eigenspace associated to
the value $h$ is one-dimensional. We shall denote by $\Phi$
the (up-to-phase) unique unit vector in this eigenspace.
The dense subspace $\fin$ consisting of the linear combinations of the
eigenvectors will be called the space of ``finite-energy'' vectors.
The representation via infinitesimal generators defines an
irreducible unitary lowest weight representation 
$\{L_n : n \in \ZZ \}$ of the Virasoro algebra satisfying
for all natural numbers $n,m$
\begin{itemize}
\item[1.] \!({\it core})
$\fin$ is a core and invariant for the closed operator $L_n$
\item[2.] \!({\it lowest weight})
if $n>0$ then $L_n \Phi = 0$
\item[3.] \!({\it unitarity})
$L_n^*=L_{-n}$
\item[4.] \!({\it Virasoro algebra relations})
on the common invariant core $\fin$
\begin{equation}
[L_n,L_m]=(n-m)L_{n+m}+\frac{c}{12}(n^3-n)\delta_{-m,n} \mathbbm 1.
\end{equation}
\end{itemize}
The correspondence between the infinitesimal generators and the
representation is the following.
For an $f\in $ Vect$(S^1) \equiv C^\infty(S^1,\RR)$ 
real vector field with Fourier coefficients
\begin{equation}
\hat{f}_n=\frac{1}{2\pi}\int e^{-in\alpha}f(e^{i\alpha})\,d\alpha 
\quad (n \in {\mathbb Z})
\end{equation}
the operator $T_0(f)$ on domain $\fin$ given by
\begin{equation}
T_0(f)=\mathop{\sum}_{n \in \mathbb Z} \hat{f}_n L_n
\end{equation} 
is well-defined and essentially self-adjoint. 
Then, denoting $T(f)$ the self-adjoint operator
obtained by the closure of $T_0(f)$ and omitting $c$ 
and $h$ indices, we have that 
\begin{equation}
e^{iT(f)} = V(\rm{Exp}(f))
\end{equation}
after an appropriate choice of the phase of the right hand side. 
Via this relation (and the conditions listed above) 
$V$ and the operators $L_n$ $(n \in \ZZ)$ completely
determine each other. $T$ is called the stress-energy tensor, it can be
looked upon as an operator valued distribution.

\smallskip

{\it Remark.}
Note that in the literature usually the operators $L_n$ $n\in\ZZ$ 
are not taken as closed operators, i.e. our notations stand for the
closure of those. 

\smallskip

The possible values of $c$ are 
$\{1-6/((m+2)(m+3))\,|\,m=1,2,3,..\}$ and $c \geq 1$, and 
for all these $h=0$ is a possible lowest weight. In case
of $h=0$, we shall denote the representation simply by 
$V_c$ (omiting the zero in the subscript), and by 
$\Omega$ the (up-to-phase) unique unit zero-energy vector 
(omiting even the subscript ``$c$''). 
For every $I \in \I$ with the
\begin{definition}
$\A_{\rm{Vir},c}(I)=\{V_{c}(\gamma)\in {\rm B}(\H_{(c,0)})| 
\,\gamma|_{I'}=\rm{id}_{I'}\}''$
\end{definition} \noindent
the net $\A_{\rm{Vir},c}$ with 
the representation of $\mob$ obtained by restriction of $V_c$ is a 
M\"obius covariant net on $S^1$, which is also diffeomorphism 
covariant with respect to the representation
$V_c$. This is the so-called Virasoro net.

With what was said before we have described all irreducible 
positive energy
representations of the diffeomorphism group: recall 
(\cite[Theorem A.1]{Carpi03b}) that
an irreducible positive energy representation of $\diff$ is equivalent
to $V_{(c,h)}$ for some value of $c$ and $h$.
The proof in \cite{Carpi03b} is based on results in \cite{loke}.

\subsection{Subnets}

A {\bf (M\"obius covariant) subnet} of the M\"obius covariant net
$(\A,U)$ is an assignment of nontrivial von Neumann algebras to the open 
proper arcs of the circle $I \mapsto \B(I)$ such that for all $I_1,I_2,I 
\in \I$ and $\varphi \in \mob$
\begin{itemize}
\item[(i)] 
$\B(I) \subset \A(I)$ 
\item[(ii)]
$I_1 \subset I_2 \Rightarrow \B(I_1) \subset \B(I_2)$ 
\item[(iii)] 
$\varphi \in \mob \Rightarrow U(\varphi)\B(I)U(\varphi)^* = 
\B(\varphi(I))$. 
\end{itemize} 
We shall use the notation $\B \subset \A$ for  subnets.

A subnet $\B \subset \A$ which is proper (namely it does not coincide 
with $\A$)   
is not a M\"obius covariant net in the precise sense of the 
definition because we do not have the cyclicity of 
the vacuum with respect to $\B$. However, this inconvenience
can be overcome by restriction to the Hilbert space 
\begin{equation}
\H_\B=\overline{\B(S^1)\Omega}=\overline{\bigvee_{I \in \I}\B(I)\Omega}
\end{equation}
where $\Omega$ is the vacuum vector.
It is evident that $\H_\B$ is invariant for $U$.
The map $I \mapsto \B(I)|_{\H_\B}$ together with the
restriction of $U$ onto $\H_\B$ is a M\"obius covariant net.
Rather direct consequences of the definition and of the properties of 
M\"obius covariant nets (such as for example the {\it Reeh-Schlieder} and 
{\it Haag property}) are: 
\begin{itemize}
\item[(i)]
for any $I \in \I$ the restriction map from $\B(I)$ to 
$\B(I)|_{\H_\B}$ is an isomorphism between von Neumann algebras,
\item[(ii)]
the map $A \mapsto P_\B A |_{\H_\B}$ where $P_\B$ is the orthogonal
projection onto $\H_\B$ and $A \in \A(I)$ for a fixed $I \in \I$ defines
a faithful normal conditional expectation from $\A(I)$ to $\B(I)$
after identifying $\B(I)$ with $\B(I)|_{\H_\B}$ using point (i),
\item[(iii)]
$\B(S^1) \cap \A(I) = \B(I)$ for all $I \in \I$.
\end{itemize}

A fundamental example of a subnet which we shall briefly describe 
here is the one determined by the stress-energy tensor in a
diffeomorphism covariant theory. Suppose the net $(\A,U)$ is
diffeomorphism covariant with respect to the strongly continuous
projective representation $V$ of $\diff$. Then the formula
\begin{equation}
\A_V(I):=\{U(\varphi): \varphi|_{I'}=\rm{id}_{I'}\}'' 
\subset \A(I)
\end{equation}
where $I \in \I$ defines a subnet $\A_V \subset \A$. 

Of course $V$ can be restricted to $\H_{\A_V}$ and the 
restriction of the subnet $\A_V$ onto this subspace
is generated by the restriction of $V$ in the sense that
local algebras are generated by the unitaries associated 
to local diffeomorphisms. Therefore by the {\it irreducibility
property} this restriction of $V$ is irreducible. Since as we have
already explained, irreducible representations are of the type 
$V_{(c,h)}$, the net $\A_V|_{\H_{\A_V}}$ is a Virasoro net
given by the representation $V|_{\H_{\A_V}}$.

Finally, let us recall an important property of the
Virasoro nets in connection with subnets. (See the
proof in \cite{Carpi98}.)
\begin{theorem}[Minimality of the Virasoro nets]
\label{minimality}
A Virasoro net does not have any proper 
M\"obius covariant subnet.
\end{theorem}
 
\section{Uniqueness in case of the Virasoro nets}
\label{sectuniqvir}
In this section we prove the uniqueness of the $\diff$-action in the case 
of Virasoro nets. This result, of itself interest, will also provide an 
important step in the uniqueness proof for the general case of 
$4$-regular nets that we shall discuss later. 

Let $\tilde{V}$ be a positive energy projective representation of
$\diff$ making $\A_{{\rm Vir},c}$ diffeomorphism covariant in the 
sense of Definition \ref{diffcov:def}. Of course in particular on the
M\"obius subgroup $\tilde{V}$ coincides with $V_c$. 
\begin{lemma} 
\label{tilde:irreducible}
$\tilde{V}$ is irreducible.
\end{lemma}
\begin{proof} 
As in the last subsection of the preliminaries, by the equation
\begin{equation}
\A_{\tilde{V}} (I)=\{\tilde{V}(\gamma): \gamma|_{I'}=\rm{id}_{I'}\}''
\;\;(I \in \I)
\end{equation}
we define a (M\"obius covariant) subnet of $\A_{{\rm Vir},c}$.
Therefore, by the minimality (cited by us in
the preliminaries as theorem \ref{minimality}) 
of the Virasoro net, taking account that
$\tilde{V}$
cannot be trivial, we have that
${\rm B}(\H_{(c,0)})=\A_{{\rm{Vir}},c}(S^1)=
\A_{\tilde{V}}(S^1)=\{\tilde{V}\}''$. 
\end{proof}

As a consequence of lemma \ref{tilde:irreducible} 
and the fact (\cite[Theorem A.1]{Carpi03b}) that the irreducible
representations are exactly the Virasoro ones, $\tilde{V}$ is an
irreducible Virasoro representation with lowest weight zero and central
charge $\tilde{c}$ (with possibly $c \neq \tilde{c}$).
We shall denote by $\{\tilde{L}_n: n\in \ZZ \}$ the 
resulting family of representing operators for the Virasoro algebra
and by $\tilde{T}$ the corresponding stress-energy tensor.

On the $\mob \subset \diff$ subgroup $\tilde{V}$ coincides with
$V_c$. Since the M\"obius vector fields are exactly the ones for which the
only nonzero Fourier coefficients are those associated to $-1,0,1$,
we have that $\tilde{L}_n=L_n$ for $n=0,1,-1$. 
The notion of ``finite-energy'' vectors (since $L_0=\tilde{L}_0$) 
is unambiguous, and any polynomial of the $L$ or $\tilde{L}$
operators is well-defined on $\fin$.
\begin{lemma}
\label{LOmega}
There exists a complex number $\zeta$ such that for every $n\in 
\ZZ$ we have $\tilde{L}_n\Omega = \zeta L_n\Omega$.
\end{lemma}
\begin{proof} 
From the theory of positive energy representations of the Virasoro algebra 
(see e.g. \cite{KaRa}) we know that every eigenvector of $L_0$ with eigenvalue 
$2$ is proportional to the nonzero vector $L_{-2} \Omega$. Since
$L_0=\tilde{L}_0$ and by the Virasoro algebra relations 
$\tilde{L}_0(\tilde{L}_{-2}\Omega)=2\tilde{L}_{-2}\Omega$, 
there must exist a
complex number $\zeta$ such that $\tilde{L}_{-2}\Omega=\zeta 
L_{-2}\Omega$.

Both vectors $\tilde{L}_{-n}\Omega$ and $L_{-n}\Omega$ vanish if 
$n<2$, so we only have to show that $\tilde{L}_{-n}\Omega=\zeta
L_{-n}\Omega$ for every integer $n\geq 2$. We do this by induction. 
For $n=2$ the equality has been shown before. 
Now assume that $\tilde{L}_{-n}\Omega=\zeta L_{-n}\Omega$
for some $n\geq 2$. Then, recalling that $L_{-1}=\tilde{L}_{-1}$ 
and using the Virasoro algebra relations we find
\begin{align*} 
(n-1)\tilde{L}_{-n-1}\Omega &= L_{-1}\tilde{L}_{-n}\Omega \\
&= \zeta L_{-1} L_{-n}\Omega \\
&= (n-1)\zeta L_{-n-1}\Omega,
\end{align*}
and the conclusion follows.
\end{proof}

We are now ready to state the main result of this section.
                                                                                       
\begin{theorem}
\label{uniqvir}
$\tilde{V}$ as projective representation coincides with $V_{c}$. 
In other words, $\A_{{\rm Vir},c}$ has a unique $\diff$ action which is
compatible with the action of $\mob$ determined by the net and its
vacuum vector $\Omega$.
\end{theorem}
\begin{proof} 
From Lemma \ref{LOmega} we find that $\tilde{T}(f)\Omega =\zeta T(f)\Omega$ 
for every real smooth function $f$. Now, if the support of $f$ is contained in 
an interval $I\in \I$ and $\Psi \in \A_{{\rm Vir},c}(I')\Omega$, it
follows 
from locality that $\tilde{T}(f)\Psi =\zeta T(f)\Psi$. But 
$\A_{{\rm Vir},c}(I')\Omega$ contains a core for $L_0$ (see the appendix
to 
\cite{Carpi99b}) and hence it is a common core for $T(f)$ and
$\tilde{T}(f)$, see e.g. \cite{BS-M} or the next section.
It follows that $\tilde{T}(f) =\zeta T(f)$ for every real smooth function $f$ 
on $S^1$ with nondense support and hence for every real
function $f$ on $S^1$. In particular, since $\tilde{L}_0=L_0$ by
assumption, we must have $\zeta=1$ and hence 
$\tilde{V}(\rm{Exp}(f))=V_{c}(\rm{Exp}(f))$ for every 
smooth real vector field $f$ on $S^1$. Our claim then follows because
$\diff$ is generated by exponentials \cite{Milnor}.
\end{proof}

\begin{corollary}
\label{isovirnets} 
Two Virasoro nets as M\"obius covariant nets are isomorphic if and only 
if they have the same central charge.
\end{corollary}

\section{Stress-energy tensor and nonsmooth vector fields}
\label{sectestimates}
Suppose we have a positive energy representation of
$\diff$. We would like to extend the representation to some
transformations that are not smooth, but still ``sufficiently regular''.
(Later we shall give more meaning to this.) For this purpose we shall
take a not necessary smooth function $f:S^1 \rightarrow \RR$ (of which we
think as a non-smooth vector field) and we will try to define a
self-adjoint operator $T(f)$ by the closure of the naive formula
$\mathop{\sum}_{n \in {\mathbb Z}} \hat{f}_n L_n$.
(As it will be discussed, even if the representation is not irreducible,
it gives rise to a corresponding representation $\{L_n: n \in {\mathbb
Z}\}$ of the Virasoro algebra with all the properties listed in the
preliminaries.)

Looking at the article of Goodman and Wallach \cite{GoWa}, we can see that
in fact everything works well with the definition of $T(f)$ even if $f$ is
not smooth but for example if $\mathop{\sum}_{n \in {\mathbb
Z}}|\hat{f}_n|(1+|n|)^3 < \infty$. Unfortunately, for the uniqueness
result we need to handle functions of less regularity. However, in the
cited article essential self-adjointness is proved by using a result in the 
paper of Nelson \cite{Ne}. Reading the work of  Nelson, we can realize
that what we really need is an $\epsilon$-independent bound on the norm
of the commutator $[\mathop{\sum}_{n \in {\mathbb Z}}\hat{f}_n L_n,
e^{-\epsilon L_0}]$ where $\epsilon > 0$. This is what we shall establish
in what follows here.

Throughout this section let $V$ be a positive energy 
representation of $\diff$. As $e^{i2\pi L_0}$ is a multiple of
the identity, the nonnegative spectrum of $L_0$ contains eigenvalues only
(at distances of integer numbers). So the linear span of the eigenvectors
$\fin$ is still a dense subspace, which we shall still call the space of
``finite-energy'' vectors. In \cite[Chapt. 1]{loke} T. Loke has shown 
that if the eigenspaces of $L_0$ are all finite dimensional than via
infinitesimal generators in the same way as it was described in the
preliminaries (but there only in the irreducible case) it gives rise to
a representation of the Virasoro algebra with a certain value of the
central charge. However, as it was pointed out in the appendix of 
\cite{Carpi03b}, this condition can be dropped. 
So the construction with infinitesimal generators works and gives us in
general a representation $\{L_n: n \in \ZZ \}$ of the Virasoro algebra
with a certain value $c>0$ of the central charge satisfying
all the properties already listed in subsection \ref{subsec:diffcov}. 
Namely, $\fin$ is an invariant core for the closed operators $L_n$ $(n \in
\ZZ)$, the adjoint operator $L_n^*$ equals to $L_{-n}$ for all $n \in \ZZ$
and on the common invariant core of the finite-energy vectors these
operators satisfy the Virasoro algebra relations.

Although equation $(2.8)$ on page no. $308$ in the article of
Goodman and Wallach \cite{GoWa} is stated for the irreducible case,
the value of the lowest weight is not involved at all and in fact
after a close look it is rather evident that as a consequence we
have in general the following:
\begin{lemma}
\label{estimate}
There exists a constant $r > 0$ independent from $k,\,v_k,n$ (but
dependent on the value of the central charge $c$) such that  
\begin{eqnarray*}
\| L_n v_k \|^2 \leq 
r^2 \left( k^2+k|n|^2+|n|^3) \right) \|v_k\|^2.
\end{eqnarray*}
where $v_k$ is an eigenvector of $L_0$ with eigenvalue
$k$ and $n \in {\mathbb Z}$. 
\end{lemma}

It is clear therefore, that $\Dl$, the domain of $L_0$, is included 
in the domain of $L_n$, and if $v \in \Dl$ then by using that
$\sqrt{1+|n|^3}\leq (1+|n|^{\frac{3}{2}})$
\begin{equation}
\label{energybound}
\|L_n v\| \leq r \,(1+|n|^{\frac{3}{2}})\,\,\|(\mathbbm 1 +L_0)v\|
\end{equation}
which is why any core for $L_0$ is a core for $L_n$. (And in particular,
as we have already stated, the finite-energy space is so.)
Related ``energy-bounds'' can be found in \cite{BS-M}. 
The above estimate has the following consequence. 
\begin{proposition}
\label{AonD}
If $a_n \in {\mathbb C} \,\, (n \in {\mathbb Z})$ is such that
$\mathop{\sum}_{n \in {\mathbb Z}}|a_n| (1+|n|^{\frac{3}{2}}) < \infty$ then 
\begin{itemize}
\item[(i)]
the operator $A=\mathop{\sum}_{n \in {\mathbb Z}}a_n L_n$ 
with domain $\Dl$ is well defined, 
(i.e. the sum strongly converges on the domain);

\item[(ii)]
if $v \in \Dl$, then as 
$N \rightarrow \infty$ the sum 
$$\mathop{\sum}_{k \in {\rm Sp}(L_0), \; k \leq N} A v_k \, 
\rightarrow A v $$ strongly, 
where the vector $v_k$ is the component of $v$ in the eigenspace of $L_0$
associated to the value $k \in $ Sp$(L_0)$,
\item[(iii)]
$A^*$ is an extension of the operator
$A^+:=\mathop{\sum}_{n \in {\mathbb Z}}\overline{a}_{-n} L_n$. 
(This again is understood as an operator with domain $\Dl$.)
\end{itemize}  
\end{proposition}
\begin{proof}
Since the sum
\begin{equation}
\label{domain}
\mathop{\sum}_{n \in {\mathbb Z}} \|a_n L_n v\|
\leq
r \left( \mathop{\sum}_{n \in {\mathbb Z}} 
 |a_n| (1+|n|^{\frac{3}{2}}) \right) \|(\mathbbm 1 +L_0)v\| < \infty 
\end{equation}
claim (i) holds. 
Claim (ii) follows from the same estimate and the fact that
$$(\mathbbm 1+L_0)\left( v - \mathop{\sum}_{k \in {\rm Sp}(L_0),
\; k \leq N} v_k \right) \rightarrow 0$$ 
as $N$ tends to $\infty$.
Finally, the last claim follows,
since for all $n$ integer $L_n^*=L_{-n}$.
\end{proof}

We now consider, for every $\epsilon >0$, the operator 
$R_{n,\epsilon}=[L_n,e^{-\epsilon L_0}]$ which is at least densely defined 
for every $n \in {\mathbb Z}$, since its domain surely contains the 
subspace $\Dl$. The following proposition gives an estimate on the norm of
this commutator which is independent of $\epsilon$. 
\begin{proposition}
\label{commutatornorm}
There exists a constant $q > 0$ independent of $\epsilon$ and $n$ such
that $\|R_{n,\epsilon}\|^2=\|[L_n,e^{-\epsilon L_0}]\|^2 \leq q |n|^3$.
\end{proposition}
\begin{proof}
For $n=0$ the statement is trivially true as $L_0$ 
commutes with any bounded function of itself. Since $L_n^*=L_{-n}$ 
and $e^{-\epsilon L_0}$ is self-adjoint it follows that
$R_{n,\epsilon} \subset -R_{-n,\epsilon}^*$, 
it suffices to demonstrate the
statement for negative values of $n$, 
and, since it also shows that $R_{n,\epsilon}$ is closable,
it is enough to verify that $\|R_{n,\epsilon}v\|^2 \leq 
q |n|^3 \|v\|^2$ whenever $v \in \fin$.

Let therefore be $n<0$, $v \in \fin$ and for every $k \in$ Sp$(L_0)$ 
let $v_k$ be again the component of the vector $v$ in the eigenspace of
$L_0$ associated to the eigenvalue $k$.
To not to get confused about
positive and negative constants, in the calculations we shall use the
positive $m:=-n$ rather than the negative $n$. Now since $L_n$ raises the
eigenvalue of $L_0$ by $m$, we have that for $k \in $Sp$(L_0)$ 
\begin{equation}
\label{Rv} 
R_{n,\epsilon}v_k=[L_n,e^{-\epsilon L_0}]v_k =
(e^{-\epsilon k}-e^{-\epsilon (k+m)})\,L_nv_k. 
\end{equation} 
The mapping $f_m: \epsilon \mapsto e^{-\epsilon k}-e^{-\epsilon (k+m)}$ 
is a positive smooth function on $\RR^+$ which goes to zero both when
$\epsilon \rightarrow 0$ and when $\epsilon \rightarrow \infty$. Therefore
$f_m$ has a maximum on $\RR^+$.
Now the only solution of the
equation $f'_m(\epsilon)=0$ is $\epsilon_m=-(1/m) \ln (k/(k+m))$. This,
together with the mentioned facts gives that 
\begin{equation} 
\label{f_m}
\mathop{\rm sup}_{\epsilon \in \RR^+}{|f_m(\epsilon)|^2}=
f_m(\epsilon_m)^2= \left(\frac{k}{k+m}\right)^{\frac{2k}{m}}
\left(\frac{m}{k+m}\right)^2 \leq \left( \frac{m}{k+m} \right)^2.
\end{equation} 
We can now return to the question of the norm of the
commutator. Equation (\ref{Rv}) shows that the vectors
$R_{n,\epsilon}v_k\,\, (k \in $Sp$(L_0))$ are in particular 
pairwise orthogonal. Using this and the fact that
only for finitely many values of $k$ the vector $v_k \neq 0$ we find
\begin{eqnarray} 
\nonumber 
\|R_{n,\epsilon}v\|^2 &=& \|R_{n,\epsilon}
\mathop{\sum}_{k \in {\rm Sp}(L_0)} v_k\|^2= 
\mathop{\sum}_{k \in {\rm Sp}(L_0)} \|R_{n,\epsilon}v_k\|^2
=\mathop{\sum}_{k \in {\rm Sp}(L_0)} |f_m(\epsilon)|^2\|L_nv_k\|^2 
\\ \nonumber
&\leq& \mathop{\sum}_{k \in {\rm Sp}(L_0)} 
\mathop{\rm sup}_{\epsilon \in \RR^+}\{|f_m(\epsilon)|^2\} \, \|L_nv_k\|^2 
\\ \nonumber 
&\leq&
\mathop{\sum}_{k \in {\rm Sp}(L_0)} \left(\frac{m}{k+m} \right)^2
r^2 (k^2+km^2+m^3) \,\|v_k\|^2 
\\ \nonumber &\leq&
\mathop{\sum}_{k \in {\rm Sp}(L_0)} r^2 (m^2+m^3+m^3) \,\|v_k\|^2 
\\
&\leq& 3 r^2 |n|^3 \,\|v\|^2, 
\end{eqnarray} 
where we have used the inequality (\ref{f_m})  and the constant $r$
is the one coming from lemma \ref{estimate} in estimating the norm square
of $L_nv_k$. 
\end{proof}

\begin{theorem} 
If $a_n \in {\mathbb C} \,\, (n \in {\mathbb Z})$ 
is such that 
$\mathop{\sum}_{n \in {\mathbb Z}}|a_n| |n|^{\frac{3}{2}} < \infty$
then $A$ is closable and
$\overline{A}=(A^+)^*$, where 
$A=\mathop{\sum}_{n \in {\mathbb Z}}a_n L_n$ and 
$A^+=\mathop{\sum}_{n \in {\mathbb Z}}\overline{a}_{-n} L_n$
considered as operators on the domain $\fin$. 
In particular, if $a_n=\overline{a_{-n}}$ for all 
$n \in {\mathbb Z}$, then $A$ is essentially self-adjoint on $\fin$. 
\end{theorem} 
\begin{proof} 
Let us first note, that because of Proposition 
\ref{AonD}, claim (iii) the operator $A$ both with domain $\Dl$ and 
$\fin$ is closable, since the domain of its adjoint surely contains
$\Dl$ and that because of proposition \ref{AonD}, claim (ii) 
$\overline{A|_{\fin}}=\overline{A|_{\Dl}}$.
Therefore from now on we shall think of $A$  
as an operator with domain $\Dl$, since in any case it 
does not change neither its closure nor its adjoint.
(Of course the same applies to the operator $A^+$.)  
Further, if $\epsilon > 0$, then the domain of the
operator $R_{A,\epsilon}=[A,e^{-\epsilon L_0}]$ is 
the whole $\Dl$ and we have that 
$R_{A,\epsilon} \subset - R_{A^+,\epsilon}^*$
where $R_{A^+,\epsilon}=[A^+,e^{-\epsilon L_0}]$. 
By using Proposition \ref{commutatornorm} with the constant 
$q$ provided by it and the condition on the sequence 
$a_n \,\, (n \in {\mathbb Z})$,
\begin{equation}
\mathop{\sum}_{n \in {\mathbb Z}}\|a_n R_{n,\epsilon}\| \leq
\mathop{\sum}_{n \in {\mathbb Z}}|a_n| q^{\frac{1}{2}} |n|^{\frac{3}{2}}
< \infty.
\end{equation}
Since $R_{A,\epsilon}=\mathop{\sum}_{n \in {\mathbb Z}}
a_n R_{n,\epsilon}$ on $\Dl$, this means that 
$\|R_{A,\epsilon}\|$ is bounded by a constant independent
of $\epsilon$. Obviously, the same is true for 
$\|R_{A^+,\epsilon}\|$. 

If $v_k$ is an eigenvector of $L_0$ with eigenvalue $k \geq 0$ 
then, as $\epsilon$ tends to zero,
\begin{equation}
R_{{A^+},\epsilon}v_k=(e^{-\epsilon k}{\mathbbm 1}-e^{-\epsilon 
L_0})A^+v_k
\rightarrow 0.
\end{equation}
Thus the operators $R_{A^+,\epsilon}$ on $\fin$ strongly converge 
to zero. Then since their norm is bounded 
by a constant independent of $\epsilon$, 
as $\epsilon \rightarrow 0$, the everywhere 
defined bounded operators 
$R_{A,\epsilon}^*=-\overline{R_{A^+,\epsilon}}$
converge strongly to zero. 

From here the proof of the theorem continues exactly 
as in \cite{Ne}, but for self-containment let us revise 
the concluding argument. Suppose $x$ is a vector in the
domain of $A^*$. Then, since $e^{-\epsilon L_0}x \in \Dl \subset 
\D (A^+)$,
we have that 
\begin{equation}
\label{eq:Nelson}
A^+ e^{-\epsilon L_0} x = A^* e^{-\epsilon L_0}x =
e^{-\epsilon L_0} A^* x  - R^*_{A,\epsilon}x.
\end{equation} 
As $\epsilon \rightarrow 0$ of course 
$e^{-\epsilon L_0}x \rightarrow x$, but now
Equation (\ref{eq:Nelson}) shows that also
$A^+ e^{-\epsilon L_0} x \rightarrow A^*x$
strongly. Therefore $A^*=\overline{A^+}$.
\end{proof}

With this we have proved the main theorem of this section. The
result ensures, that if the continuous function 
$f: S^1 \rightarrow \RR$ with Fourier coefficients
$\hat{f}_n$ $(n \in \ZZ)$ is such, that the norm 
\begin{equation} 
\|f\|_{\frac{3}{2}}=
\mathop{\sum}_{n \in {\mathbb Z}}|\hat{f}_n|(1+|n|^{\frac{3}{2}})
\end{equation}
is finite, 
then $\mathop{\sum}_{n \in {\mathbb Z}} \hat{f}_n L_n$ is an essentially
self-adjoint operator on $\fin$. As in the case of smooth functions,
we will denote by $T(f)$ the corresponding self-adjoint operator 
obtained by taking closure. We continue by investigating the
continuity property of the stress-energy tensor $T$.
\begin{proposition}
\label{convergence} 
For every continuous real function $f$ on $S^1$ of finite
$\|\cdot\|_{\frac{3}{2}}$ norm and for every $v \in \Dl$
we have
\begin{equation}
\|T(f)v\| \leq r \|f\|_{\frac{3}{2}} \|(\mathbbm 1+L_0)v\|
\end{equation}
where $r$ is the positive constant appearing in 
Lemma \ref{estimate}. Moreover, if $f$ and $f_n \,\, (n\in \NN)$ 
are continuous real functions on $S^1$ of finite
$\|\cdot\|_{\frac{3}{2}}$ norm, and $\|f_n-f\|_{\frac{3}{2}}$ 
converges to zero as $n$ tends to $\infty$, then $T(f_n) 
\rightarrow T(f)$ in the strong resolvent sense. In
particular, $e^{iT(f_n)} \rightarrow e^{iT(f)}$ strongly.
\end{proposition}
\begin{proof}
The claimed inequality is an immediate consequence of 
the inequality in Eq. (\ref{energybound}) and the
definition of the $\|\cdot\|_{\frac{3}{2}}$ norm.
Now by this estimate for every $v \in \Dl$ we have
that $T(f_n)v$ converges to $T(f)v$. Since $\Dl$
is a common core for these self-adjoint operators,
the conclusion follows (see e.g. \cite[Sect. VIII.7]{RSI}).
\end{proof}
In the next section we shall need to determine the geometrical
properties of the adjoint action of $e^{iT(f)}$ for a certain $f$
nonsmooth vector field. If $f$ was smooth, we would know what 
the unitary $e^{iT(f)}$ ``does'' since it is the operator
associated by the representation to the diffeomorphism Exp$(f)$.
Thanks to the last proposition, to obtain information in the
case when $f$ is not smooth, all we will have to do is
to approximate it with smooth ones in an appropriate way.
As it will be clear later, the following lemma shows that 
for our purposes the smooth vector fields are ``many enough''.
\begin{lemma}
\label{density}
Let $I \subset S^1, \,\, I\neq\emptyset$ be an open 
interval (or even the whole circle). If
$f$ is a real continuous function of finite 
$\|\cdot\|_{\frac{3}{2}}$ norm with support 
contained in $I$, then there exists a sequence
$f_k,\,(k=1,2,..)$ of real smooth functions with 
support still in $I$ such that 
$\lim_{k\to \infty}\|f_k-f\|_{\frac{3}{2}}=0$.
\end{lemma}
\begin{proof}
The proof follows standard arguments relying on 
convolution with smooth functions. 
Let $\varphi_k \,\, (k=1,2,..)$ be a sequence
of positive smooth functions on $S^1$ with support
shrinking to the point $1 \in S^1$ such that for
all $k \in \NN$ 
$$\frac{1}{2\pi}\int_0^{2\pi} \varphi_k(e^{i\alpha})
\,d\alpha=1.$$
Then for $k$ large enough the convolution $\varphi_k*f$
is a smooth real function with support in $I$. Moreover we have 
$$\|\varphi_k*f-f\|_{\frac{3}{2}} \leq \sum_{n\in
\ZZ}|(\hat{(\varphi_k)}_n-1) \hat{f}_n|\,(1+|n|^{\frac{3}{2}})$$  
where the left-hand side goes to zero when $k \to \infty$ since
$|(\hat{\varphi_k})_n| \leq 1$ and 
$\lim_{k\to\infty}(\hat{\varphi_k})_n=1$.
\end{proof}
\section{Uniqueness in case of $4$-regularity}
\label{sectuniq4reg}
Suppose $(\A,U)$ is a M\"obius covariant net on the circle which is also
diffeomorphism covariant with the representation $V$ of $\diff$.
Further, suppose that there exists another positive energy 
representation $\tilde{V}$ of $\diff$  which 
also makes $(\A,U)$ diffeomorphism covariant in the sense 
of Definition \ref{diffcov:def}.
The representations $V$ and $\tilde{V}$, as before, via infinitesimal
generators, give rise the representations $L$ and $\tilde{L}$ of the
Virasoro algebra (with possibly different values of the central charge).
The corresponding stress-energy tensor fields we shall denote by $T$ and
$\tilde{T}$. As the two projective representations coincide on the
M\"obius subgroup, $L_n=\tilde{L}_n$ for $n=-1,0,1$ or to put it in
another way, $T(g)=\tilde{T}(g)$ whenever $g$ is a M\"obius vector field.

Considering the two representations we can define two M\"obius covariant
subnets: the subnet $\A_V$ defined for every open proper arc 
$I \subset S^1$ by
\begin{equation}
\A_V(I)=\{V(\varphi)\,|\,\varphi|_{I'}=\rm{id}_{I'}\} \subset \A(I),
\end{equation}
and the subnet $\A_{\tilde{V}}$ defined similarly but with the
representation $V$ replaced by $\tilde{V}$. The restriction of the subnet
$\A_V$ onto the closed linear subspace $\overline{\A_V\Omega}$, as it has
been cited from \cite{Carpi03b} several times by now, is a Virasoro net 
for a certain value of central charge $c$ and so it will be denoted by
$\A_{\rm{Vir}}$.  In the same way the restriction of $\A_{\tilde{V}}$ onto 
$\overline{\A_{\tilde{V}}\Omega}$ is another Virasoro net 
(with the possibly different value of central charge
$\tilde{c}$) and will be denoted by $\A_{\tilde{\rm Vir}}$. But a
Virasoro net is a {\it minimal} net, and this gives a very
strong restriction on the possible ways the two subnets $\A_V$ and
$\A_{\tilde{V}}$ can ``differ''.
\begin{proposition}
\label{localintersection1}
If $\A_V$ and $\A_{\tilde{V}}$ as subnets are not equal, then
for all $I \subset S^1$ open proper arc
$\A_V(I) \cap \A_{\tilde{V}}(I) = \CC {\mathbbm 1}$.
\end{proposition}
\begin{proof} Suppose that $\A_V(I) \cap \A_{\tilde{V}}(I)$ is nontrivial 
for  a given (and hence, by M\"obius covariance, for all) $I\in \I$. 
Then the subnet $I \mapsto \A_V(I) \cap \A_{\tilde{V}}(I) \subset \A(I)$ 
when restricted to $\overline{\A_V\Omega}$ is a M\"obius covariant 
subnet of $\A_{\rm{Vir}}$, therefore by minimality (cited by us as Theorem
\ref{minimality}) it must coincide with $\A_{\rm Vir}$. 
On the other hand, for an open proper arc $I$ the restriction map from
$\A_V(I)$ to $\A_{\rm{Vir}}(I)$ is an isomorphism. 
So we have that $\A_V(I) \cap \A_{\tilde{V}}(I)$ coincides with
$\A_V(I)$ for every $I\in \I$. 
But of course by interchanging $V$ and $\tilde{V}$, it must also coincide 
with $\A_{\tilde{V}}(I)$  for every $I\in \I$ and this concludes the 
proof.
\end{proof}
\begin{proposition}
\label{localintersection2}
If $\A_V$ and $\A_{\tilde{V}}$ as subnets are equal then so are 
$V$ and $\tilde{V}$ as projective representation; i.e. 
$\rm{Ad}(V(\varphi))=\rm{Ad}(\tilde{V}(\varphi))$
for all $\varphi \in \diff$.
\end{proposition}
\begin{proof}
By the condition of the proposition the representation 
$\tilde{V}$ can be restricted onto
$\overline{\A_V\Omega}$ and this gives
a positive energy representation
of $\diff$, which is compatible
with the Virasoro net $\A_{\rm{Vir}}$ (as M\"obius covariant net).  
Hence, by the uniqueness result for the Virasoro nets (Theorem
\ref{uniqvir}) it must be equal (as a projective representation) 
with the restriction of $V$ onto the same subspace. But if $\varphi$ 
is a diffeomorphism ``localized'' in the open proper arc
$I \subset S^1$, that is, $\varphi|_{I'}=\rm{id}_{I'}$,
then --- since $V(\varphi),\tilde{V}(\varphi) \in \A_V(I)$,
and the restriction map from $\A_V(I)$ to $\A_{\rm{Vir}}(I)$
is an isomorphism --- the operator $V(\varphi)^*\tilde{V}(\varphi)$ must
be a multiple of the identity. This is enough for the equality, since
$\diff$ is generated by localized diffeomorphisms.
\end{proof}

Thus, by the previous two propositions if the local intersections
of the subnets $\A_V$ and $\A_{\tilde{V}}$ are not trivial, then we have
that $V$ and $\tilde{V}$, as projective representations, are equal.
We know that the intersection of the two {\it algebras} $\A_V(S^1)$ and
$\A_{\tilde{V}}(S^1)$ cannot be trivial: it contains the unitaries
associated to M\"obius transformations, for example. 
Unfortunately, there is no M\"obius transformation --- apart from the
identity --- that would be local. However, we can construct local
transformations that are {\it piecewise M\"obius}. 
Naturally, they will not be smooth, but, as we shall see it, 
by choosing the parameters rightly, we can achieve 
once differentiability, with discontinuities (``jumps'') appearing 
at the endpoints of the pieces only in the second derivative.
We have essentially three things to do:
\begin{itemize}    
\item[(i)]
we have to construct such a $\zeta$ piecewise M\"obius transformation,

\item[(i)] 
we must show, that although $\zeta$ is not smooth, it is
sufficiently regular so that the expressions $V(\zeta)$ and
$\tilde{V}(\zeta)$ are meaningful, and finally,

\item[(ii)]
we must show that the adjoint action of these unitaries on the 
algebras corresponding to some pieces is completely determined 
(since the geometrical part of this action on each piece is
M\"obius), and we must investigate that under what condition
on the net it implies that the two unitaries are in fact multiples
of each other. 
\end{itemize}
We begin with the construction of a piecewise M\"obius transformation.
For a $z \in S^1$ let $I_{(z,iz)} \subset S^1$ be the open quarter-arc
with
endpoints $z$ and $iz$. The real M\"obius vector field $g_1$ given by
the formula
\begin{equation}
g_1(z)=(i-1)z+2-(i+1)z^{-1}
\end{equation}
is zero in the two points $1,i \in S^1$. Hence by setting 
$g_p \,\,(p=1,i,-1,-i)$ to be the real M\"obius vector field
determined by the equation
\begin{equation}
g_p(pz)=p^2 g_1(z)
\end{equation}
the map
\begin{equation}
z \mapsto \left\{ 
\begin{matrix}
g_1(z) & {\rm if}\,\,  z \in I_{(1,i)} \\ 
g_i(z) & {\rm if}\,\,  z \in I_{(i,-1)} \\ 
g_{-1}(z) & {\rm if}\,\, z \in I_{(-1,-i)} \\ 
g_{-i}(z) & {\rm if}\,\, z \in I_{(-i,1)} \\ 
\end{matrix} \right.
\end{equation}
defines a unique continuous function $f: S^1 \rightarrow \RR$.
We shall think of this function as a nonsmooth vector field. 
It has four points at which it is zero: the points $1,i,-1$
and $-i$. On each of the four quarter-arc between these points it 
coincides with a M\"obius vector field (of course on each arc
with a different one). 
\begin{lemma}
$\|f\|_{\frac{3}{2}} < \infty$
\end{lemma}
\begin{proof}
By direct calculation not only $f$, but also its derivative is
continuous. Its second derivative is of course still smooth on
each of the four open intervals, and at the endpoints it has only 
finite ``jumps''. Therefore it is a function of bounded variation, 
and hence there is a constant $M>0$ such that the absolut value of its 
Fourier coefficient associated
to any integer $n$ is bounded by
$\frac{M}{|n|}$ (see \cite[Sect. I.4]{katznelson})
which in turn implies that
$|\hat{f}_n| \leq \frac{M}{|n|^3}$ for all $n \in \ZZ$. 
\end{proof}
This means that we can consider the self-adjoint operators 
$T(f)$ and $\tilde{T}(f)$. By construction, $T(f)$ is affiliated to 
$\A_V(S^1)$ and $\tilde{T}(f)$ is affiliated to $\A_{\tilde{V}}(S^1)$.
\begin{proposition}
For all $t \in \RR$ the adjoint actions of $e^{itT(f)}$ and 
$e^{it\tilde{T}(f)}$ restricted to the algebra $\A(I_{(p,ip)})$
coincide with that of $e^{itT(g_p)}=e^{it\tilde{T}(g_p)}$,
where $p=1,i,-1,-i$. 
\end{proposition}
\begin{proof}
The continuous real function $f-g_p$ is zero on $I_{(p,ip)}$.
Since its $\|\cdot\|_{\frac{3}{2}}$ norm is finite, we can consider
the operator $T(f-g_p)$ which then by Lemma \ref{density} and Proposition
\ref{convergence} is affiliated to $\A(I_{(p,ip)}')$. On the common
core of the finite-energy vectors 
\begin{equation}
T(f)=T(f-g_p)+T(g_p).
\end{equation}
Then, since $T(f-g_p)$ is affiliated to the commutant of $\A(I_{(p,ip)})$ 
and
\begin{equation}
\rm{Ad}(e^{itT(g_p)})(\A(I_{(p,ip)})=\A(I_{(p,ip)})
\end{equation}
for all $t \in \RR$, by a simple use of the Trotter product-formula 
(e.g. \cite[Theorem VIII.31]{RSI}) we 
have the part of the proposition concerning the adjoint action of
$e^{itT(f)}$. Similar argument justifies the assertions for
$e^{it\tilde{T}(f)}$. 
\end{proof}
This means that on the algebra associated to the 
four open quarter-arc the adjoint actions of
$e^{itT(f)}$ and $e^{it\tilde{T}(f)}$ coincide. 
Hence if $\A$ is at least $4$-regular, the unitary
$e^{itT(f)}e^{-it\tilde{T}(f)}$ must be a multiple of the identity. 
It follows that $e^{itT(f)}e^{-itT(g_1)} \in \A_V(S^1)$
and $e^{it\tilde{T}(f)}e^{-it\tilde{T}(g_1)} \in \A_{\tilde{V}}(S^1)$ 
are multiples of each other and --- since they act trivially on
$\A(I_{(1,i)})$ --- that they belong to the local intersection
$\A_V(I_{(1,i)}') \cap \A_{\tilde{V}}(I_{(1,i)}')$ for all real $t$. 
But of course they cannot be just multiples of the identity: for example
because $T(f)\Omega$ is a nonzero vector (the real $f$ is not a
M\"obius vector field, so it cannot have zero Fourier coefficients
associated to all values of $n < -1$) which is orthogonal to $\Omega$.
Then, by Prop. \ref{localintersection1} and Prop.
\ref{localintersection2} we can conclude that $V$ and $\tilde{V}$, as
projective representations are equal.

Thus we have proved that
\begin{theorem}
\label{uniq4-reg} 
Let $(\A,U)$ be an at least $4$-regular diffeomorphism
covariant net on the circle. Then there is a unique 
projective representation $V$ of $\diff$  
which makes $(\A,U)$ diffeomorphism covariant
in the sense of Definition \ref{diffcov:def}.
\end{theorem}

Let us formulate now some important consequences of the fact that the
whole representation $V$ must already be encoded in the M\"obius covariant
net (with its given representation of the M\"obius group, or
equivalently,
with its given vacuum vector). Remember that a positive energy
representation of $\diff$ always gives rise to a representation of the 
Virasoro algebra
(see the discussion in the beginning of Section \ref{sectestimates}), so
in particular it always has a central charge.
\begin{corollary}
\label{invariants}
Let $(\A,U)$ be a $4$-regular net with the representation
$V$ of $\diff$ making it diffeomorphism covariant. Then the
representation class of $V$, and in particular its central charge
$c>0$ is an invariant of the M\"obius covariant net $(\A,U)$.
\end{corollary}

Another interesting thing to note here is the model-independent 
proof for the commutation between internal symmetries and diffeomorphism
symmetry.
\begin{definition}
A unitary $W$ on the Hilbert space $\H_\A$ is called an (unbroken) 
{\bf internal symmetry} of the net $(\A,U)$ if for every $I \in \I$
\begin{equation}
W \A(I) W^*=\A(I)
\end{equation}
and $W\Omega=\Omega$ where $\Omega$ is the vacuum vector of $(\A,U)$.
\end{definition}
By our uniqueness theorem we can state the following conclusion.
\begin{corollary}
\label{com_int} 
Let $W$ be an internal symmetry of the net $(\A,U)$
having diffeomorphism symmetry. If $\A$ is at least $4$-regular,
than the unique representation $V$ of $\diff$ making the
net diffeomorphism covariant must commute with $W$. 
\end{corollary}
\begin{proof}
Since $W$ commutes with the representation $U$ (see \cite{FrG})
the projective representation $WVW^*$ of $\diff$ still makes 
the net $(\A,U)$ diffeomorphism covariant. Hence by Theorem
\ref{uniq4-reg} it must coincide with $V$. It follows that
for every $\gamma \in \diff$ the unitary 
$WV(\gamma)W^*V(\gamma)^*$ is a multiple
$\lambda(\gamma)$ of the identity, and in fact 
it turns out that the complex valued function
$\gamma \mapsto \lambda(\gamma)$ is a character of 
the group $\diff$. But the latter is a simple 
noncommutative group (see e.g. \cite{Milnor}), 
and hence $\lambda$ is trivial.
\end{proof}
\section{Infinite tensor products and nets admitting no
diffeomorphism symmetry}
\label{secttensorproduct}
In this section we shall use our uniqueness results to exhibit a class of 
M\"obius covariant nets on $S^1$ that definitely do not have a
diffeomorphism symmetry. 

Let $(\A_n, U_n)$, $n=1,2,...$ be a sequence of M\"obius covariant nets on 
$\s1$ and let $\Omega_n$, $n=1,2,...$ be the corresponding sequence of 
vacuum vectors. We can define the infinite tensor product net 
\begin{equation}
\A \equiv \bigotimes_n \A_n 
\end{equation}
on the 
(separable) infinite tensor product 
Hilbert space 
$$\H_\A:=\bigotimes_n^{(\Omega_n)}\H_{\A_n}$$ 
by 
\begin{equation}
\A(I):= \bigotimes_n\A_n(I),
\end{equation}
cf. \cite{VN38}.
It is fairly easy to show that $\A$ together with the representation 
$\bigotimes_n U_n$ is a M\"{o}bius covariant net on $S^1$ which is strongly 
additive (resp. n-regular) when each net $\A_n$, $n=1,2 ...$ is strongly 
additive (resp. n-regular). 

We shall need the following proposition which is of interest of its 
own.

\begin{proposition} 
\label{difftensor}
Let, $\A, \B$, two 4-regular M\"obius 
covariant nets on $\s1$. 
If $\A$ and $\A \otimes \B$ are diffeomorphism covariant and 
$V_\A, V_{\A\otimes \B}$ are the corresponding representations of
$\diff$, then $\B$ is diffeomorphism covariant with a representation 
$V_\B$ satisfying $V_\A \otimes V_\B =  V_{\A\otimes \B}.$
\end{proposition}
\begin{proof}
Let us consider the M\"obius covariant net $\A\otimes\A \otimes\B$. 
By assumption it is a 4-regular diffeomorphism covariant net on $\s1$
and the corresponding representation of $\diff$ is given by 
$V:=V_\A \otimes  V_{\A\otimes \B}$. Let $F$ be the unitary operator 
on $\H_\A\otimes \H_\A \otimes \H_\B$ which flips the first two components 
of the tensor product. It is easy to see that $F$ is an internal symmetry 
of the net $\A\otimes\A \otimes\B$ and hence,
by Corollary \ref{com_int}, it must commute with $V$. Since, 
for every $C \in {\rm B}(\H_\A)$, $\gamma \in \diff$,  we have 
$$V(\gamma)\left(C\otimes {\mathbbm 1} \otimes {\mathbbm 1}
\right)V(\gamma)^* 
=\left(V_\A(\gamma)CV_\A(\gamma)^*\right)\otimes {\mathbbm 1} 
\otimes {\mathbbm 1},$$ we find 
\begin{align*}
V(\gamma)\left({\mathbbm 1}\otimes C \otimes {\mathbbm 1}
\right)V(\gamma)^* &=
FV(\gamma)F^*\left({\mathbbm 1}
\otimes C \otimes {\mathbbm 1}
\right)FV(\gamma)^*F^* \\
&= FV(\gamma)\left(C\otimes {\mathbbm 1}
\otimes {\mathbbm 1}
\right)V(\gamma)^*F^* \\
&= F \left( \left(V_\A(\gamma)CV_\A(\gamma)^*\right)\otimes 
{\mathbbm 1} \otimes {\mathbbm 1}
\right)F^* \\
&= {\mathbbm 1} \otimes\left(V_\A(\gamma)CV_\A(\gamma)^*\right)
\otimes {\mathbbm 1}.
\end{align*}

It follows that, for every $\gamma \in \diff$, 
$$V(\gamma)\left(V_\A(\gamma)\otimes V_\A(\gamma) 
\otimes {\mathbbm 1} \right)^*
\in \left({\rm B}(\H_\A) \otimes {\rm B}(\H_\A)  
\otimes {\mathbbm 1} \right)'= 
{\mathbbm 1} \otimes {\mathbbm 1}
\otimes {\rm B}(\H_\B)$$
and hence that there is a projective unitary representation $V_\B$ of 
$\diff$ on $\H_\B$ such that 
$$V=V_\A \otimes V_\A \otimes V_\B.$$
Then, the conclusion easily follows. 
 
\end{proof}

We are now ready to state the main result of this section. 

\begin{theorem}
\label{nondifftensor}
Let $\A_n$, $n\in \NN$ be a sequence of 4-regular 
diffeomorphism covariant nets on $\s1$. Then the infinite tensor product
net $\otimes_n \A_n$ together with the corresponding tensor product
representation of {\rm $\mob$} is not diffeomorphism covariant.
\end{theorem}
\begin{proof} We denote by $V_n$ the representation of $\diff$ corresponding 
to $\A_n$, $n\in \NN$, and by $c_n$ its central charge. For every 
positive integer $k$ the net $\otimes_n \A_n$ is isomorphic to 
$$\left(\A_1\otimes ...\otimes \A_k \right)\otimes \B_k,$$
where
$$\B_k := \bigotimes_{n>k}\A_k.$$
Let us assume that $\otimes_n\A_n$ is diffeomorphism covariant and let
$V$ be the corresponding representation of $\diff$. By Prop. 
\ref{difftensor} there is a positive energy 
representation $V_{\B_k}$ of $\diff$ making $\B_k$ diffeomorphism covariant 
and such that 
$$V =V_1\otimes...\otimes V_k \otimes V_{\B_k}.$$
Hence the central charge $c$ of $V$ satisfies 
$$c=c_1+...+c_k+ c(\B_\k) \geq \frac{k}{2},$$
where $c(\B_\k)$ is the central charge of $V_{\B_\k}$, since
the minimal possible value for a central charge is $1/2$. 
However, by the arbitrariness of $k$ we have a contradiction and 
the conclusion follows.
\end{proof}
The examples of non diffeomorphism covariant nets considered in 
\cite{koester03a} are not strongly additive. However, they satisfy the 
trace class condition, namely $e^{-\beta L_0}$ is a trace class operator for 
every $\beta >0$ and hence they have the split property by
\cite[Theorem 3.2]{D'ALR}. Conversely one can use Theorem 
\ref{nondifftensor} to give many examples of non diffeomorphism covariant 
strongly additive M\"obius covariant nets on $\s1$. In these examples 
the operator $e^{-\beta L_0}$ is not compact for every value of $\beta$ 
since the eigenvalue 2 of $L_0$ always appears with infinite multiplicity.
In fact, e.g. if the sequence $c_k$ of the central charges contains a 
constant subsequence the infinite tensor product net 
$\otimes_n \A_n$ does not satisfy the split property as a consequence of 
\cite[Theorem 9.2]{WSSI}. 

\bigskip

\noindent{\bf Acknowledgements.} We would like to thank Roberto Longo for 
suggesting the problem on the uniqueness of the diffeomorphism symmetry
and for useful discussions. Some of the results contained in this paper 
were announced by the second named author (M. W.) at the conference 
on ``Operator Algebras and Mathematical Physics'', held in Sinaia
(Romania) in July, 2003. He would also like to thank the organizers for
the invitation.

\end{document}